\magnification1200
\overfullrule=0pt
\def\O{{\cal O}}

\def\o{\omega}
\def\P{{\cal P}}
\def\Pr{{\rm Pr}}
\def\R{{\cal R}}
\def\A{{\cal A}}

\def\C{{\cal C}}
\def\qed{{\it qed}}
\def\mapdown#1{\downarrow\rlap{\hbox{$#1$}}}
\null
\vskip2cm

\centerline{\bf ON THE HILBERT SCHEME OF  A PRYM VARIETY}

\medskip
\centerline{\bf H. Lange \ \   \  \ E. Sernesi}
\vskip1cm
\hfill {\it Dedicated to the memory of  Fabio Bardelli}

\vskip1cm

{\bf Introduction} 

\bigskip

We work over the field of complex numbers. In this paper we consider 
the Prym map $\P: \R_g \to \A_{g-1}$ from the moduli space of unramified 
double covers of projective irreducible and nonsingular curves of genus 
$g \ge 6$ to the moduli space of principally polarized abelian varieties
 of dimension $g-1$. 
If $\pi:\tilde C \to C$  is such a double cover with $C$ non hyperelliptic we 
consider  the natural embedding $\tilde C \subset P$ (defined up to 
translation) of $\tilde C$ into the Prym variety $P$ of $\pi$ and we study 
the local structure of the Hilbert scheme $Hilb^P$ of $P$ at the point 
 $[\tilde C]$ (here and through the paper we adopt the notation $[-]$ for the 
point of a 
moduli space or of a Hilbert scheme which parametrizes the object $-$).
We show that this structure is related with the local geometry of the Prym 
map, 
or more precisely with the validity of the  infinitesimal version of Torelli's 
theorem for Pryms at 
$[\pi]$ (see \S 3 for the definitions). 

The  results we prove are the following.

\bigskip

{\bf   Proposition} {\sl If the infinitesimal Torelli theorem for Pryms
holds at 
$[\pi]$ then $Hilb^P$ is nonsingular of dimension $g-1$ at $[\tilde C]$
(i.e. $\tilde C$ is unobstructed)
and the only  deformations of $\tilde C$ in $P$ are translations.}

\bigskip
It is known that if the Clifford index of $C$ is at least 3 then the 
condition of the Proposition is satisfied. Therefore we have in particular:

\bigskip

{\bf   Corollary} {\sl If Cliff$(C) \ge 3$ then 
$Hilb^P$ is nonsingular of dimension $g-1$ at $[\tilde C]$
(i.e. $\tilde C$ is unobstructed)
and the only  deformations of $\tilde C$ in $P$ are translations.}

\bigskip
On the other side we have the following result:

\bigskip

{\bf   Theorem} {\sl Assume that the following conditions are satisfied:

\noindent
(a) the infinitesimal Torelli theorem for Pryms fails  at $[\pi]$;

\noindent
(b) $[\pi]$ is an isolated point of the fibre $\P^{-1}(\P([\pi]))$.

Then  $\tilde C$ is obstructed. Moreover   the only local deformations of 
$\tilde C$ in $P$ are translations 
and the only irreducible component of $Hilb^P$ containing $[\tilde C]$ is 
everywhere non reduced. 

Conversely,  if $\tilde C$ is obstructed then the infinitesimal Torelli 
theorem for 
Pryms fails at $[\pi]$.}

\bigskip

Using  these results we give some examples in which 
$\tilde C \subset P$ is obstructed and  some  
in which we have unobstructedness but the infinitesimal Torelli theorem for 
Pryms fails. 
The examples we construct are obtained 
from double covers belonging to $\R_6$ and to $\R_7$. For their construction
we make use of a result proved in \S 4 which is a slight extension of a 
theorem 
of Recillas (see [R]).

\bigskip

The paper is divided into 5 sections. In \S 1 we discuss the Hilbert scheme  
of curves Abel-Jacobi embedded in their jacobian. We prove  that such 
curves are obstructed precisely when they are hyperelliptic of genus
$g \ge 3$.  This case is not relevant for what follows but it is worth keeping 
in mind the analogies between the two cases.
In \S 2 we consider our 
problem and we study the  conditions for the unobstructedness of $[\tilde C]$. 
We use the well 
known cohomological description of  certain tangent spaces and of maps 
between them. In \S 3 we relate these results with the infinitesimal Torelli 
theorem and 
we prove our main result. In \S 4 we give a proof of the   extension of 
Recillas' theorem. The 
final \S 5 contains the examples. 

\bigskip

\vskip1cm
{\bf 1. The case of curves in their jacobians} 

\bigskip
 Consider a projective nonsingular 
irreducible curve $C$
of genus $g \ge 2$, let $JC := Pic^0(C)$ be the jacobian variety of $C$, and 
let $j: C \to JC$ be an Abel-Jacobi map. We want to study the local structure 
of the Hilbert scheme $Hilb^{JC}$ of $JC$ at $[j(C)]$ (the point 
parametrizing $j(C)$).
Since $j$ is an embedding we will identify $C$  with $j(C)$. We have an exact 
sequence of locally free sheaves on $C$:
$$
0 \to T_C \to T_{JC|C} \to N_C \to 0 \leqno (1)
$$
We have a canonical isomorphism 
$T_{JC|C} \cong H^1(\O_C)\otimes \O_C$
and therefore the cohomology sequence of (1) is:
$$
0 \to H^1(\O_C)  \to H^0(N_C) 
\buildrel \delta \over \longrightarrow H^1(T_C)
\buildrel \sigma \over \longrightarrow H^1(\O_C)\otimes H^1(\O_C)
\to H^1(N_C) \leqno (2)
$$
The family of translations of $C$ in $JC$ is parametrized by $JC$ itself, and 
the map $H^1(\O_C)  \to H^0(N_C)$ in (2) is precisely the characteristic map 
of this family at the point $0$. Therefore we have the following Lemma, whose 
proof is obvious:

\bigskip

{\bf 1.1 Lemma} {\sl The following conditions are equivalent:

\noindent
(a) $Hilb^{JC}$ is nonsingular of dimension $g$ at $[C]$.

\noindent
(b) $\delta = 0$

\noindent
(c) $\sigma$ is injective

\noindent
If these conditions are satisfied the only local deformations of $C$ in $JC$
are  translations.}

\bigskip
Using this Lemma we can prove the following

\bigskip

{\bf 1.2 Theorem } {\sl Suppose that $C$ has genus $g \ge 3$.  

\noindent
 (a) If $C$ is non hyperelliptic then $Hilb^{JC}$
is nonsingular of dimension $g$ at $[C]$.

\noindent
(b) If $C$ is hyperelliptic then the connected component of $Hilb^{JC}$ 
containing $[C]$ is irreducible of dimension $g$ and everywhere non reduced  
with Zariski 
tangent space of dimension $2g-2$.

In both cases the only deformations of $C$ in $JC$ are translations.}

  \medskip
  {\bf Proof.}  The transpose of $\sigma$ is the multiplication map:
  $$
  \sigma^\vee: H^0(\o_C)\otimes H^0(\o_C) \to H^0(\o_C^{\otimes 2})
  $$
   (see [G], Lemma 3). This map, by Noether's theorem, is surjective if $C$ is 
non 
hyperelliptic and 
has corank $g-2$ if $C$ is hyperelliptic (see [ACGH]). Therefore, in view of 
Lemma 1.1, part (a) follows.
   
   Now assume that $C$ is hyperelliptic and that $\bar C \subset JC$ is a 
closed subscheme such that $[\bar C]$ belongs to the connected component of 
  $Hilb^{JC}$ containing $[C]$. By  the criterion of Matsusaka-Ran (see [LB])
  $\bar C = C_1 \cup \cdots \cup C_r$ is a reduced curve of compact type, and 
$JC$ and $JC_1\times\cdots\times JC_r$ are isomorphic as ppav's. Then it 
follows that $r=1$ and $\bar C$ is irreducible and nonsingular because $C$ 
is.     
   Now we apply Torelli's theorem to conclude that $\bar C$ is a translate of 
$C$. It follows that the connected component of $Hilb^{JC}$ containing $[C]$ 
is irreducible of dimension $g$ and parametrizes the translates of $C$. On 
the other hand by (2) we have $h^0(N_C)= 2g-2 > g$. The conclusion follows. 
\qed
   
  \bigskip
  Theorem 1.2 can be interpreted in terms of the Torelli morphism
  $$
  \tau: M_g \to A_g
  $$
   from the moduli stack of projective nonsingular curves of genus $g$ to the 
moduli stack of principally polarized abelian varieties of dimension $g$.
   The surjectivity of $\sigma^\vee$ is equivalent to that of the 
multiplication map
$$
 S^2H^0(\o_C) \to H^0(\o_C^{\otimes 2})
$$
which is the codifferential of $\tau$ at $[C]$. Hence the surjectivity of this 
map is equivalent to the infinitesimal Torelli theorem for $C$ (see [OS]). 
Therefore Theorem 1.2 implies the following:

\bigskip
{\bf 1.3 Corollary } {\sl $C$ is unobstructed in $JC$ if and only if the 
infinitesimal Torelli theorem holds for $C$.}

\bigskip
{\sl Remarks.}  (i) 
The proof of Theorem 1.2(a) appeared already in [G], but the argument does not 
appear to be complete. A proof is 
also given in [Bl] using the semiregularity map, but it is more 
complicated; moreover the 
semiregularity map does not seem to be able to detect what happens in case 
(b).

(ii)  In the case $g=2$ we have that $C$ is unobstructed in $JC$ because the 
semiregularity map $H^1(N_C) \to H^2(\O_{JC})$ is injective since 
$H^1(\O_{JC}(C))=0$ by the ampleness of $C$ in $JC$.

\vskip1cm

{\bf 2. The Hilbert scheme of the Prym variety at $\bf [\tilde C]$}

\bigskip

Let now $\pi : \tilde C \to C$ be an unramified double cover of a projective 
nonsingular irreducible curve $C$ of genus $g \ge 3$, so that $\tilde C$ has 
genus $\tilde g = 2g-1$.  Let $\eta \in Pic^0(C)$ be the  2-division point 
corresponding to $\pi$. We have a canonical isogeny
$J\tilde C \to JC\times P$, where $P$ is the Prym variety of $\pi$.
{\it Throughout this section we assume $C$ to be non hyperelliptic}. 
Under this hypothesis we have an embedding $\alpha: \tilde C \to P$ which is 
obtained as the composition
$$
\tilde C \to J\tilde C \to JC\times P \to P
$$
(see [LB]). We will identify $\tilde C$ with $\alpha(\tilde C)$. We want to 
study the Hilbert scheme $Hilb^P$ locally at the point $[\tilde C]$. 

In analogy with the situation studied in \S 1, we consider the exact sequence
of locally free sheaves on $\tilde C$:
$$
0 \to T_{\tilde C} \to T_{P|\tilde C} \to N_{\tilde C} \to 0 \leqno (3)
$$
We have a canonical isomorphism
$T_{P|\tilde C} \cong H^1(C,\eta)\otimes \O_{\tilde C}$
so that the cohomology sequence of (3) becomes:
$$
0 \to H^1(C,\eta) \to H^0(\tilde C,N_{\tilde C})
\buildrel \delta \over \longrightarrow H^1(\tilde C,T_{\tilde C})
\buildrel \sigma \over \longrightarrow 
H^1(C,\eta)\otimes H^1(\tilde C,\O_{\tilde C})
\to H^1(\tilde C,N_{\tilde C}) \leqno (4)
$$

Along the same lines of \S 1 we can state the following

\bigskip

{\bf 2.1 Lemma} {\sl 
The following conditions are equivalent:

\noindent
(a) $Hilb^P$ is nonsingular of dimension $g-1$ at $[\tilde C]$.

\noindent
(b) $\delta = 0$

\noindent
(c) $\sigma$ is injective

\noindent
If these conditions are satisfied the only local deformations of $\tilde C$ in
$P$  are translations.}

\bigskip
In order to understand  the conditions of  Lemma 2.1 
 we must study the map $\sigma$, or equivalently its transpose
$\sigma^\vee$. We have canonical isomorphisms:
$$ H^1(\tilde C,\O_{\tilde C})^\vee \cong H^0(\tilde C,\o_{\tilde C})
\cong H^0(C,\o_C)\oplus H^0(C,\o_C\otimes\eta)
$$
and
$$
H^1(\tilde C,T_{\tilde C})^\vee \cong H^0(\tilde C,\o_{\tilde C}^{\otimes 2})
\cong 
H^0(C,\o_C^{\otimes 2}\otimes \eta)\oplus H^0(C,\o_C^{\otimes 2})
$$
corresponding to the decompositions into $+1$ and $-1$ eigenvalues under the 
action induced by the involution on $\tilde C$. Hence
$$
\sigma^\vee: H^0(C,\o_C\otimes \eta)\bigotimes 
[H^0(C,\o_C)\oplus H^0(C,\o_C\otimes\eta)] \to
H^0(C,\o_C^{\otimes 2}\otimes \eta)\oplus H^0(C,\o_C^{\otimes 2})
$$
and it is induced by multiplication of sections ([B], page 382). Therefore, 
after decomposing 
the domain of $\sigma^\vee$ as 
$$
H^0(C,\o_C\otimes \eta)\bigotimes 
[H^0(C,\o_C)\oplus H^0(C,\o_C\otimes\eta)] = $$
$$ =
[H^0(C,\o_C\otimes \eta)\otimes H^0(C,\o_C)] \bigoplus
[H^0(C,\o_C\otimes \eta)\otimes H^0(C,\o_C\otimes\eta)]
$$
we see that  $\sigma^\vee = \mu \oplus \nu$ where:
$$
\mu: H^0(C,\o_C\otimes \eta)\otimes H^0(C,\o_C) \to 
H^0(C,\o_C^{\otimes 2}\otimes \eta)
$$
and
$$
\nu: H^0(C,\o_C\otimes \eta)\otimes H^0(C,\o_C\otimes\eta) \to
H^0(C,\o_C^{\otimes 2})
$$

The following Lemma is well known (see [B], Prop. 7.7):

\bigskip
{\bf 2.2 Lemma} {\sl The sheaf $\o_C\otimes \eta$ is not very ample if and 
only 
if there exist points $x,y,z,t \in C$ such that $\eta\cong \O_C(x+y-z-t)$. If 
these conditions are satisfied then the map $\nu$ is not surjective.}

\bigskip
{\bf 2.3 Proposition} {\sl  (i) In each of the following cases the map $\mu$ 
is 
surjective: \hfill\break
(a)  $C$ is not bielliptic
\hfill\break
(b)  $\nu$ is surjective.
\hfill\break
(ii) If $C$ is not bielliptic then
${\rm {\rm cork}}(\sigma^\vee) = {\rm cork}(\nu)$. 
In particular if $C$ is not bielliptic the 
surjectivity of $\sigma^\vee$ is equivalent to the surjectivity of $\nu$.}

\medskip
{\bf Proof.} (i) Note first that in both cases (a) and (b) the linear series  
$|\o_C\otimes \eta|$
is base point free and is not composed with an involution: in fact,  
in case (a) since $C$ is not hyperelliptic $|\o_C\otimes \eta|$ is base point 
free; moreover  if it were composed with an involution then, since  
deg$(\o_C\otimes \eta)=2g-2$, the morphism 
$$\phi_\eta: C \to {\bf P}^{g-2}
$$
would be of degree 2 onto a curve of degree $g-1$, which has genus $\le 1$, a 
contradiction. In case (b)  the assertion is true by Lemma 2.2.

 Let ${\underline b} := P_1+\cdots +P_{g-3} \in C^{(g-3)}$ be 
general. Consider the exact sequence on $C$:
$$
0 \to \o_C\otimes\eta(-{\underline b}) \to \o_C\otimes \eta \to {\cal T}
\to 0
$$
where ${\cal T}$ is a torsion sheaf supported on ${\underline b}$. 
Multiplying 
firstly by $H^0(\o_C)$ and taking cohomology, and secondly by $\o_C$ and 
taking cohomology, we obtain the following commutative diagram with exact 
rows,  where the vertical maps are given by multiplication:
$$
\matrix{
0\to&H^0(\o_C\otimes\eta (-{\underline b}))\otimes H^0(\o_C)&\to&
H^0(\o_C\otimes\eta)\otimes H^0(\o_C)&\to&
H^0({\cal T})\otimes H^0(\o_C)&\to 0 \cr
&\mapdown{\mu_{{\underline b}}}&&\mapdown{\mu}&&\mapdown{\bar\mu}&\cr
0\to&H^0(\o_C^{\otimes 2}\otimes\eta (-{\underline b}))&\to&
 H^0(\o_C^{\otimes 2}\otimes\eta)&\to& H^0({\cal T}\otimes \o_C)&\to 0}
 $$
Since $|\o_C\otimes \eta|$ is not composed with an involution
and ${\underline b}$ is generic, 
$\o_C\otimes\eta (-{\underline b})$ is base point free, and
by the base point free pencil trick we find:
$$
{\rm ker}( \mu_{{\underline b}} ) = H^0({\underline b}\otimes \eta) = 0
$$
hence:
$${\rm rk}(\mu_{{\underline b}}) = 2g = 
h^0(\o_C^{\otimes 2}\otimes\eta (-{\underline b}))
$$
i.e.  $\mu_{{\underline b}}$ is surjective. On the other hand $\bar \mu$ is 
surjective because $\o_C$ is globally generated. The conclusion follows from 
the above diagram. 

(ii)  follows immediately from part (i) and from the relation between the maps 
$\sigma, \mu, \nu$. \qed

\bigskip

Collecting all we have said so far we can state the following:

\bigskip

{\bf 2.4 Corollary} {\sl If $\nu$ is surjective  then 
$Hilb^P$ is nonsingular of dimension $g-1$ at $[\tilde C]$
(i.e. $\tilde C$ is unobstructed)
and the only local deformations of $\tilde C$ in $P$ are translations.}

\bigskip

As an application  we can prove the following: 

\bigskip

{\bf 2.5 Corollary} {\sl If Cliff$(C) \ge 3$ then 
$Hilb^P$ is nonsingular of dimension $g-1$ at $[\tilde C]$
(i.e. $\tilde C$ is unobstructed)
and the only local deformations of $\tilde C$ in $P$ are translations.}

\medskip

{\bf Proof.}   It  follows easily from a result of   [GL] (see e.g. [LS]) 
that if 
Cliff$(C) \ge 3$ then  the map $\nu$ is surjective.
Therefore the Corollary follows from Corollary 2.4.
\qed

\vskip1cm

{\bf 3. $\bf Hilb^P$ and the infinitesimal Torelli theorem
for  Pryms}

\bigskip

 We keep the notations of \S 2. Consider the Prym morphism:
$$
\P: \R_g \to \A_{g-1}
$$
which goes from the coarse moduli scheme of \'etale double covers of curves of 
genus 
$g$ to the coarse moduli scheme of ppav   of dimension $g-1$, $g \ge 6$.  
These schemes have singularities due to the presence of automorphisms of the 
objects they classify. Therefore if we want to study the infinitesimal 
properties of $\P$ it is more natural to consider the corresponding 
moduli stacks 
 $R_g, A_{g-1}$.
 The  Prym construction defines a morphism of stacks
 $$
 \Pr: R_g \to A_{g-1}
 $$
 Then the map $\nu$ considered in \S 2 coincides with the codifferential of 
 Pr at $[\pi]$ (see [B], Prop. 7.5, which implies this statement modulo 
obvious modifications).  
Therefore the surjectivity of $\nu$ is equivalent   to $\Pr$ being an 
immersion at $[\pi]$ (see [B], 7.6). In this case
 we say that 
the  infinitesimal Torelli theorem for Pryms holds at $[\pi]$, according to 
the terminology most commonly used nowadays.
In view of Corollary 2.4 we can therefore state the following:

\bigskip

{\bf 3.1 Proposition} {\sl If the infinitesimal Torelli theorem for Pryms 
holds at $[\pi]$ 
then $Hilb^P$ is nonsingular of dimension $g-1$ at $[\tilde C]$
(i.e. $\tilde C$ is unobstructed)
and the only local deformations of $\tilde C$ in $P$ are translations.}

\bigskip

In the case 
Cliff$(C) \le 2$ the infinitesimal Torelli theorem  for Pryms in general 
fails, i.e. in 
general $\nu$ is not surjective. Our next goal is to relate the 
obstructedness of 
$\tilde C$ in $P$ to the failure of the infinitesimal Torelli theorem for 
Pryms.
The main result of this section  is the following:

\bigskip

{\bf 3.2 Theorem} {\sl Assume that the following conditions are satisfied:

\noindent
(a) the infinitesimal Torelli theorem for Pryms fails at $[\pi]$;

\noindent
(b) $[\pi]$ is an isolated point of the fibre $\P^{-1}(\P([\pi]))$.

Then  $\tilde C$ is obstructed. Moreover   the only local deformations of 
$\tilde C$ in $P$ are translations; 
in particular the only irreducible component of $Hilb^P$ containing $[\tilde 
C]$ 
is 
everywhere non reduced of dimension $g-1$. 

Conversely,  if $\tilde C$ is obstructed then the infinitesimal Torelli 
theorem fails at 
$[\pi]$.}

\medskip
{\bf Proof.}
 By (a) the map $\delta$ in the exact sequence (4) is non zero. Assume by 
contradiction that $[\tilde C]$ is unobstructed. Then we can 
find a nonsingular curve $S \subset Hilb^P$ passing through $[\tilde C]$ such 
that 
   $\delta(T_{S,[\tilde C]}) \ne 0$. This condition implies that the 
functorial 
morphism $S \to {\cal M}_{\tilde g}$ 
    defined by the family of curves 
   ${\cal C} \to S$ (which can be assumed to be smooth)  
 is not constant and therefore   this family does not consist of curves all 
isomorphic to $\tilde C$. But this is impossible: in fact for each  curve 
$\tilde C'$ in the 
family we have 
$$
\tilde C'  \equiv_{{\rm num}} \tilde C 
\equiv_{{\rm num}} {2\over (g-2)!}\Xi^{g-2}
$$
 so that by a 
theorem of Welters (see [W]) there is an \'etale double cover 
$\pi':\tilde C' \to C'$, $(P,\Xi)$ is the 
Prym variety of $\pi'$ and $\tilde C'$ is Prym embedded. But this 
contradicts condition (b) if $\tilde C'$ is not isomorphic to $\tilde C$ 
because $[\pi'] \in \P^{-1}(\P([\pi]))$.

This analysis also shows that locally the only  deformations of $\tilde C$ in 
$P$ are translations; and since $\delta\ne 0$ the Zariski  tangent space of 
$ Hilb^P$ at $[\tilde C]$ has dimension larger than $g-1$. This proves also 
the 
last assertion.

The converse is a special case of Prop. 3.1.
\qed

\bigskip
The Theorem does not say anything in the case when $[\pi]$ is a non
 isolated point of the fibre $\P^{-1}(\P([\pi]))$. We will see in \S 5 that in 
this case there are examples where $[\tilde C]$ is unobstructed.

\vskip1cm

{\bf 4. Further considerations.}

\bigskip
A Theorem of Recillas [R] says that if $\pi: \tilde C\to C$ is a double cover 
with  $C$ trigonal (but not 
hyperelliptic) of genus $g$ then  $\P([\pi])=[JX]$ with $X$  a 4-gonal 
curve, 
and the pair $(X,g^1_4)$ is uniquely determined. 
A consequence of this result and of a theorem of Mumford [M] which gives a 
list of the Prym varieties which are jacobians, is that 
$$
\P^{-1}([JX])=W^1_4(X) \leqno (5)
$$ 
set-theoretically if $g-1 \ge 6$. 
This  says in particular that for $g \ge 11$ the Prym map is 1-1 
on 
$\R_{g,T}$ (=the locus of \'etale double covers of trigonal curves): this 
follows 
from the fact that $W^1_4(X)$ consists of at most one point  if $g(X)\ge 10$.
If $g-1=5$ then (5) is not true but we have a strict inclusion $\supset$ (see 
\S 4).
The following 
Proposition gives some further information which will suffice 
for some applications. 


\bigskip

 {\bf 4.1 Proposition} {\sl    Assume that $X$ is a nonsingular irreducible 
curve of genus $g-1 \ge 5$, non hyperelliptic nor trigonal,  and   such that 
every $g^1_4$ on $X$ has no divisors of the form $2P+2Q$ or $4P$ and let
$\pi: \tilde C \to C$ be an unramified double cover, with $C$ trigonal of 
genus $g$, such that $\P([\pi])= [JX]$. Then there is a canonical  
isomorphism 
between the kernel of the differential of $\Pr$ at $[\pi]$ and the Zariski 
tangent space of $W^1_4(X)$ at the line bundle $L$ corresponding to $\pi$.}

\medskip
{\bf Proof.} 
Since $X$ is not trigonal    we may view 
$W^1_4(X)$ as parametrizing   4-1 morphisms of $X$ into 
${\bf P}^1$.  Let $\varphi: X \to {\bf P}^1$ be the 4-1 cover defined by $L$.
 
 Let $B= Spec({\bf C}[\epsilon])$ and consider a family of deformations of 
$\varphi$ parametrized by $B$:
$$
\matrix{
 X \times B& &\to &&   {\bf P}^1\times B\cr
&\searrow && \swarrow &\cr
&& B& &&} 
$$
To this family we can associate a family of deformations of $\pi$ just 
extending Recillas' construction, as follows. Consider the second relative 
symmetric product over 
 ${\bf P}^1\times B$:
 $$
 \tilde {\cal C} :=S^{(2)}_{{\bf P}^1\times B}(X \times B)
 $$
 which comes endowed with an induced family of morphisms of degree 6:
 $$ 
 \matrix{
 f^{(2)}:&& \tilde {\cal C} &&\to &&{\bf P}^1\times B \cr
 &&&\searrow && \swarrow &\cr
&&&& B& &}
 $$
 On $\tilde {\cal C}$ there is a natural involution $\iota$ commuting with 
the projection to $B$. Letting ${\cal C} = \tilde {\cal C}/\iota$, we obtain 
a family  parametrized by $B$:
 $$
 \matrix{
 \tilde {\cal C} && \buildrel \Pi \over \longrightarrow && {\cal C} \cr
 &\searrow && \swarrow &\cr
 & & B&& } \leqno (6)
 $$
 such that $f^{(2)}$ factors through $\Pi$. Therefore (6) is a first order 
deformation of $\pi$ ; moreover (6) is contained in 
$\P^{-1}([JX])$  by construction and therefore it is an element of 
${\rm ker}(d\Pr_{[\pi]})$.
 
 Conversely, assume a family (6) given,  and assume that (6) is contained in 
${\rm ker}(d\Pr_{[\pi]})$. Then we  also have a family of triple 
covers of
 ${\bf P}^1$:
 $$
 \matrix{
  \C &&\to &&{\bf P}^1\times B\cr
 &\searrow && \swarrow &\cr
 && B&& }
 $$ 
    
 Correspondingly we have an inclusion 
 ${\bf P}^1\times B \subset S^{(3)}_B(\C)$, and an 
\'etale morphism of degree 8  
 $$
 \Pi^{(3)}: S^{(3)}_B(\tilde\C) \to S^{(3)}_B(\C)
 $$
 induced by $\Pi$. Let ${\cal D} :=\Pi^{(3)-1}({\bf P}^1\times B)$. The 
involution  $\iota$ on 
 $\tilde \C$ induces an involution on $S^{(3)}_B(\tilde\C)$ which commutes 
with  $\Pi^{(3)}$ and induces an involution on ${\cal D}$. We obtain a 
commutative diagram:
 $$
 \matrix{
 {\cal D}& \to & {\bf P}^1\times B \cr
 \downarrow & \nearrow & \downarrow \cr
 {\cal D}/\iota& \to & B } \leqno (7)
 $$
 where the diagonal morphism defines a family of deformations of $X$ 
  with an assigned $g^1_4$ on the family. The assumption that (6) is 
contained in ${\rm ker}(d\Pr_{[\pi]})$ means
that the family of ppav obtained as Pryms of the family (6) is the family of 
jacobians of ${\cal D}/\iota \to B$ and that it is trivial. By the 
infinitesimal Torelli  theorem  for jacobians
 ${\cal D}/\iota \to B$ is the trivial family as well. Therefore diagram (7) 
gives us a family of deformations of $\varphi$.  \qed

  \vskip1cm

{\bf 5. Examples}

\bigskip

{\bf 5.1 }  Let $\bar X \subset {\bf P}^2$ be an irreducible sextic having 4 
distinct nodes $N_1,\ldots,N_4$, and let $X$ be the normalization of 
$\bar X$, which has genus 6. If no three among $N_1,\ldots,N_4$ are on a 
line then $W^1_4(X)$ consists of five nonsingular points, the $g^1_4$'s cut 
by the four pencils of lines through each of the nodes and by the pencil of 
conics containing $N_1,\ldots,N_4$. From Proposition 4.1 it then follows that
  the infinitesimal Torelli theorem holds at the five double covers
  $\pi: \tilde C \to C$ of  trigonal curves   of genus 7 such that 
  $\P([\pi])=[JX]$. 
  
  Assume now that $\bar X$ has three of the nodes, say $N_1,N_2,N_3$, on a 
line. 
Then the $g^1_4$ defined by $N_4$ and that defined by the pencil of conics 
are identified to a unique element $L$ of $W^1_4(X)$ with a 1-dimensional 
Zariski tangent space. Applying Proposition 4.1 to this case we see that 
the infinitesimal Torelli theorem fails at the double cover $\pi$ 
corresponding to 
$(X,L)$  in the fibre 
  $\P^{-1}([JX])$. 
  Moreover, since equality (5) implies that $\P^{-1}([JX])$ is finite, from 
Theorem 3.2 we deduce that $\tilde C$ is obstructed in $JX$ and that the only 
component of $Hilb^P$ containing $[\tilde C]$ is everywhere non reduced and 
consists of translates of $[\tilde C]$.
  
  A count of parameters shows that in this way we get a 14-dimensional locus 
where the infinitesimal Torelli theorem fails inside the 18-dimensional space 
$\R_7$.

  \bigskip
  
  {\bf 5.2}  Let's consider the Prym map $\P: \R_6 \to \A_5$. This case has 
been extensively studied in [DS] and offers a wide variety of examples, but 
it is  not yet completely understood from the point of view of the 
infinitesimal Torelli 
theorem.  
  Recall that both domain and codomain are irreducible of dimension 15.
  Some loci where the infinitesimal Torelli theorem fails are the following.
  
  \medskip
  5.2.1 - Consider a nonsingular curve $C\subset {\bf P}^4$ obtained as a 
general 
hyperplane section of a Reye congruence in ${\bf P}^5$, i.e. of an Enriques 
surface $S$ of degree 10 contained in a nonsingular quadric. Then $C$ is a 
curve of genus 6, embedded with a Prym canonical linear series
  $|\o\otimes\eta|$;  since $C$ is contained in a quadric it follows that the 
map $\nu$ is not surjective, and therefore the infinitesimal Torelli theorem 
fails at 
the double cover 
$\pi:\tilde C \to C$ associated to $\eta$. 
  
  Naranjo-Verra  proved that the fibre 
  $\P^{-1}(\P([\pi]))$ is discrete [NV]. Therefore from Theorem 3.2 it follows 
that $Hilb^P$ is obstructed at $\tilde C$.
  
  Note that a  count of 
parameters shows that
  the locus of double covers $\pi$ constructed in this way has dimension
  $14 = 9+5$ (9 for the moduli of Enriques surfaces and 5 for the hyperplane 
sections), i.e. it is a divisor in $\R_6$. In particular it follows that a 
general such curve $C$ is not trigonal since trigonal curves depend on 13 
parameters.
  
  \medskip
  
  5.2.2 - Consider an irreducible sextic $\bar C\subset {\bf P}^2$ with four 
nodes 
such 
that two of its bitangents meet in one of the nodes, say $N$. Then the 
normalization $C$ has 
genus 6 and the $g^1_4$ defined by the pencil of lines through $N$ has two 
divisors of the form $2P+2Q$. It follows that there is a 2-division point 
$\eta \in Pic(C)$ such that $\o\otimes \eta$ is not very ample 
and the map $\nu$ is not surjective (use Lemma 2.2).
Therefore the infinitesimal Torelli theorem fails at the double cover 
$\pi:\tilde C \to 
C$ associated 
to $\eta$.  
  The locus in $\R_6$ defined by this family of examples is disjoint from the 
previous one because there the line bundles $\o\otimes \eta$ were very 
ample.  It is not clear to us what the dimensions of the fibres
  $\P^{-1}(\P([\pi]))$ are and therefore whether $[\tilde C]$ is obstructed in 
this case.
  
  \medskip
  
  5.2.3 - Another locus where the infinitesimal Torelli fails is $\R_{6,T}$, 
the locus of 
double covers 
of trigonal curves. It has dimension 13, and the restriction of $\P$ to  
$\R_{6,T}$ has general   fibre of dimension 1, as it follows from Recillas' 
Theorem
recalling that $W^1_4(X)$ for a curve $X$ of genus 5 has dimension 1.  
  
  What is interesting here is that $W^1_4(X)= \Theta_{sing}$, the singular 
locus of the theta divisor of $JX$:  for a  general $X$ this is a 
nonsingular curve of genus 11 which has an involution $\iota$ with quotient a 
nonsingular plane quintic $C$. The  
  double cover $\pi: W^1_4(X) \to C$ is associated to a 2-division point 
$\eta$ such that $\O(1)\otimes\eta$ is an even theta-characteristic (see [DS] 
for details). Moreover $\P([\pi])=[JX]$ again, by [M]. Therefore we see that 
for 
a general $X$ of genus 5 we have:
  $$
  \P^{-1}([JX]) = W^1_4(X) \cup \{[\pi]\}
  $$
  In particular  the fibre of $\P$ is not equidimensional. Moreover $\nu$ is 
surjective at $[\pi]$  (see [DS], part II, \S 5) and therefore the curve 
  $W^1_4(X)= \Theta_{sing}$ is unobstructed in $JX$. Note that this gives an 
example of a double cover $\pi$ of a curve of Clifford index 1 (namely a 
nonsingular plane quintic) at which the infinitesimal Torelli theorem holds. 

Note also that, since ${\rm cork}(\nu)$ is 1-dimensional 
 if $\pi: \tilde C \to C$ is a double cover of a general trigonal curve of 
genus 
6, 
 we have that ${\rm cork}(\sigma)$ is 1-dimensional as well, by Prop. 2.3 
(clearly a general trigonal $C$ is not bielliptic). Therefore 
$\delta$ has rank 1 and
 $$
 h^0(\tilde C,N_{\tilde C}) = g
 $$
 With some 
extra effort one can easily show that in this case $\tilde C$ is unobstructed 
in $JX$.  In fact consider a small (1-dimensional) neighborhood $A$ of $[\pi]$ 
in the fibre 
$\P^{-1}([JX])$ and let 
$$
\matrix{
\tilde \C & \subset & JX \times A \cr
\downarrow && \cr
A &&} \leqno (8)
$$
be the corresponding 1-parameter family of deformations of $\tilde C$ in 
$JX$. Since this family has varying moduli, in the exact sequence (4) we have 
$0 \ne\delta(v)\in H^1(\tilde C,T_{\tilde C})$ if $v\ne 0$ is a tangent 
vector to $A$ at the point $a_0$ parametrizing $\tilde C$, and 
$\delta(v)$ generates Im$(\delta)$. Now consider a 
small 
  neighborhood $B$ of $0$ in $JX$ and build a new family:
 $$
\matrix{
\tilde \C' & \subset & JX \times A\times B \cr
\downarrow && \cr
A\times B &&}
$$ 
 whose fibre over $(a,b)$ is the curve $t_b^ \ast(\tilde \C_a)$, i.e. the 
translate by $b$ of the fibre  $\tilde \C_a$ of the family (8).  It is clear 
that the characteristic map 
 $$
 T_{A\times B,(a_0,0)} \to H^0(\tilde C, N_{\tilde C})
 $$ 
 is an isomorphism, proving that $\tilde C$ is unobstructed. Note that 
 $h^0(\tilde C, N_{\tilde C}) > g-1$ in this case, and $\tilde C$ has non 
trivial moduli.
 
\vskip1cm

\centerline{\bf References}

\bigskip

[ACGH] Arbarello et al.: {\it Geometry of Algebraic Curves I}, Springer
Grundlehren 267 (1985).

[B] Beauville A.: Varietes de Prym et Jacobiennes intermediares, {\it 
Ann. Sci. Ecole Norm. Sup.}, 10 (1977), 309-391.

[Bl] Bloch S.: Semi-regularity and De Rham cohomology, {\it Inventiones 
Math.} 17 (1972), 51-66.

[DS] Donagi R. - Smith R.C.: The structure of the Prym map, {\it Acta Math.} 
146 (1982), 25-102.

[G] Griffiths Ph.: Some remarks and examples on continuous systems and 
moduli, {\it J. of Math. and Mech.} 16 (1967), 789-802.

[GL] Green M. - Lazarsfeld R.: On the projective normality of complete linear 
series on an algebraic curve, {\it Inventiones Math.} 83 (1986),
73-90.

[LB] Lange H.- Birkenhake Ch.: {\it Complex Abelian Varieties}, Springer 
Grundlehren 
302 (1992).

[LS]  Lange H. - Sernesi E.: Quadrics containing a Prym-canonical curve.
{\it J. Algebraic Geometry}, 5 (1996), 387-399.

[M] Mumford, D.: Prym varieties I, in {\it Contributions to Analysis}, N.Y  
Academic Press 1974.

[NV] Naranjo J.C.- Verra A.: In preparation.

[OS] Oort F. - Steenbrink J.: On the local Torelli problem for algebraic 
curves, {\it Jour. Geom. Alg. Angers 1979}, Sijhoff and Noordhoff (1980), 
157-204.

[R] Recillas S.: Jacobians of curves with a $g^1_4$ are Prym varieties of 
trigonal curves, {\it Bol. Soc. Mat. Mexicana} 19 (1974), 9-13.

[W] Welters G.: Curves of twice the principal class on principally polarized 
abelian varieties. {\it Indag. Math.} 94 (1987), 87-109.

\bigskip\noindent
Mathematisches Institut
 
 \noindent
Bismarckstr. $1{1\over 2}$, D-91054 Erlangen (Germany)

\noindent
lange@mi.uni-erlangen.de

\bigskip\noindent
Dipartimento di Matematica, Universit\`a Roma Tre

\noindent
L.go S.L. Murialdo 1, 00146 Roma (Italy)

\noindent
sernesi@mat.uniroma3.it

\end